\documentclass[5p,twocolumn]{elsarticle}
\journal{Automatica}
\usepackage[utf8]{inputenc}
\usepackage{amssymb}
\usepackage{graphicx}
\usepackage{color}

\usepackage{capt-of}
\usepackage{eqnarray,amsmath}
\usepackage{mathtools}
\usepackage{enumerate} 
\usepackage{pgfplots}
\usepackage{longtable}
\usepackage{tikz}

\newtheorem{assumption}{Assumption}
\newtheorem{lemma}{Lemma}
\newtheorem{mydef}[lemma]{Definition}
\newtheorem{theorem}[lemma]{Theorem}
\newtheorem{corollary}[lemma]{Corollary}
\newtheorem{remark}{Remark}
\newtheorem{proposition}[lemma]{Proposition}
\newtheorem{problem}{Problem}
\newtheorem{conjecture}[lemma]{Conjecture}
\usepackage{multirow}
\usepackage{subcaption}
\usepackage{standalone}

\usepackage{fontawesome}

\newlength{\figureheight}
\newlength{\figurewidth}
\definecolor{colorLeader}{rgb}{0.00000,0.44700,0.74100}%
\definecolor{colorFollower}{rgb}{0.85000,0.32500,0.09800}%

\newcommand{\legendLeader}{\raisebox{2pt}{\tikz{\draw[color=colorLeader, line width=1.5pt] (0,0) -- (5mm,0);}}}
\newcommand{\legendFollower}{\raisebox{2pt}{\tikz{\draw[color=colorFollower, line width=1.5pt] (0,0) -- (5mm,0);}}}
\newcommand{\legendFollowerSimulated}{\raisebox{2pt}{\tikz{\draw[color=darkgray, line width=1.5pt] (0,0) -- (5mm,0);}}}

\let\leq\leqslant
\let\geq\geqslant

\newcommand{\calA}{\ensuremath{\mathcal{A}}}
\newcommand{\calB}{\ensuremath{\mathcal{B}}}
\newcommand{\calC}{\ensuremath{\mathcal{C}}}

\newcommand{\calE}{\ensuremath{\mathcal{E}}}

\newcommand{\calI}{\ensuremath{\mathcal{I}}}

\newcommand{\calK}{\ensuremath{\mathcal{K}}}
\newcommand{\calL}{\ensuremath{\mathcal{L}}}

\newcommand{\bmat}{\begin{matrix}}
\newcommand{\emat}{\end{matrix}}
\newcommand{\bbm}{\begin{bmatrix}}
\newcommand{\ebm}{\end{bmatrix}}
\newcommand{\bbma}{\begin{bmatrix*}[r]}
\newcommand{\ebma}{\end{bmatrix*}}
\newcommand{\bpm}{\begin{pmatrix}}
\newcommand{\epm}{\end{pmatrix}}
\newcommand{\bpma}{\begin{pmatrix*}[r]}
\newcommand{\epma}{\end{pmatrix*}}
\newcommand{\bvm}{\begin{vmatrix}}
\newcommand{\evm}{\end{vmatrix}}
\newcommand{\bse}{\begin{subequations}}
\newcommand{\ese}{\end{subequations}}
\newcommand{\beq}{\begin{equation}}
\newcommand{\eeq}{\end{equation}}
\newcommand{\ben}{\renewcommand{\labelenumi}{\arabic{enumi}.}
\renewcommand{\theenumi}{\arabic{enumi}}\begin{enumerate}}

\newcommand{\een}{\end{enumerate}}

\newcommand{\beni}{\renewcommand{\labelenumi}{\roman{enumi}.}
\renewcommand{\theenumi}{\roman{enumi}}\begin{enumerate}}

\newcommand{\eeni}{\end{enumerate}}

\newcommand{\bena}{\renewcommand{\labelenumi}{\alph{enumi}.}
\renewcommand{\theenumi}{\alph{enumi}}\begin{enumerate}}

\newcommand{\eena}{\end{enumerate}}

\newcommand{\bit}{\begin{itemize}}
\newcommand{\eit}{\end{itemize}}
\newcommand{\bthe}{\begin{theorem}}
\newcommand{\ethe}{\end{theorem}}
\newcommand{\blem}{\begin{lemma}}
\newcommand{\elem}{\end{lemma}}
\newcommand{\bprop}{\begin{proposition}}
\newcommand{\eprop}{\end{proposition}}
\newcommand{\bex}{\begin{example}}
\newcommand{\eex}{\end{example}}
\newcommand{\bas}{\begin{assumption}}
\newcommand{\eas}{\end{assumption}}
\newcommand{\bre}{\begin{remark}}
\newcommand{\ere}{\end{remark}}
\newcommand{\bcor}{\begin{corollary}}
\newcommand{\ecor}{\end{corollary}}
\newcommand{\bdfn}{\begin{definition}}
\newcommand{\edfn}{\end{definition}}
\newcommand{\bcon}{\begin{conjecture}}
\newcommand{\econ}{\end{conjecture}}

\newcommand{\R}{\ensuremath{\mathbb R}}
\newcommand{\C}{\ensuremath{\mathbb{C}}}

\newcommand{\BP}{\noindent{\bf Proof. }}
\newcommand{\EP}{\hspace*{\fill} $\blacksquare$\bigskip\noindent}

\newcommand{\figref}[1]{\figurename~\ref{#1}}
\newcommand{\tabref}[1]{\tablename~\ref{#1}}

\newcommand{\mb}{\begin{bmatrix}}
	\newcommand{\mbb}{\end{bmatrix}}

\newcommand{\mbs}{\begin{smallbmatrix}}
	\newcommand{\mbbs}{\end{smallbmatrix}}

\newenvironment{smallbmatrix}
{\left[\begin{smallmatrix}}
	{\end{smallmatrix}\right]}
\usepackage{xcolor}

\usepackage{lipsum}

\pgfplotsset{compat = 1.18}

\newcommand\blfootnote[1]{
    \begingroup
    \renewcommand\thefootnote{}\footnote{#1}
    \addtocounter{footnote}{-1}
    \endgroup
}

\begin{document}
	\begin{frontmatter}
		
		\title{Analysis and experimental validation of decentralized controllers for delayed spacing policies in vehicle-platooning}

		\author{Paul Wijnbergen 
        \corref{mycorrespondingauthor}}
		\cortext[mycorrespondingauthor]{Corresponding author}
		\ead{p.wijnbergen@tue.nl}

        \author{Redmer de Haan}
        \author{Erjen Lefeber}
        
		\address{The authors are with the department of Mechanical Engineering, Eindhoven University of Technology,
Eindhoven, The Netherlands }

		\begin{abstract}
				In this paper, a novel approach to spacing policies for vehicle platoons and a framework for control design is presented. Whereas traditional approaches aim to mitigate the effect of actuation and communication delays on the spacing error, this paper presents spacing policies that account for delays. A framework for decentralized control design is presented and necessary and sufficient conditions for tracking of delayed-based spacing policies are stated. It is shown that due to the control design, string stability is induced by the spacing policy. The results are supported by experimental validations of the theoretical results.
		\end{abstract}
		
		
	\end{frontmatter}

\section{Introduction}
Vehicle platooning amounts to the formation of closely-spaced groups of vehicles and has the potential to increase road safety, improve traffic flow, and reduce fuel consumption, see the overviews \cite{varaiya_1993,alam_2015b,besselink_2016,wang_2020}. Consequently, the automatic control of vehicles for platooning has been studied extensively, \textit{e.g.}, \cite{levine_1966,stankovic_2000,fax_2004,zhang_2016,LefeberPloegNijmeijer20}.
\blfootnote{This work was supported by the
European Union’s Horizon Europe program under grant agreement No
101069748 – SELFY project.}

A crucial aspect of control strategies for vehicle platooning is defining the desired distance between consecutive vehicles, commonly referred to as the \emph{spacing policy}. Well-known examples are the constant spacing policy \cite{swaroop_1999}, where a constant distance between vehicles is desired, and the constant headway spacing policy \cite{ioannou_1993,swaroop_1994}, where the desired spacing is dependent on the velocity of the follower vehicle. 

Having the ability to maintain short distances in a platoon has potential benefits such as increased road capacity \cite{Xiao_2018} and reduced fuel consumption \cite{Alam_2010}. However, delays have shown to be detrimental to platoon performance, with string stability frequently failing to manifest in practice \cite{liu2001effects, xiao2008stability, devika2019control}. Therefore, not taking these delays into account in the controller design can compromise the associated benefits.

To deal with delays in platoons, there are approaches that use an approximation of the delay in controller design \cite{Xing_2016, HaanSandLef24}, or approaches that use a predictor based approach to compensate the actuation delay \cite{Davis_2021, Bekiaris_2023}. However, in these approaches, the actuation delays are primarily treated as impediments. Fully mitigating the effects requires state predictions of the preceding vehicle, which pose practical challenges due to communication delays.

In this paper, we take a novel approach by integrating actuator delays not only into the vehicle model but also into the spacing policies. We present a framework that allows for straightforward control design and enables us to deal with actuation delays. This perspective simplifies the stability analysis and, more importantly, aligns experimental results with theoretical predictions.

The contributions of this paper are the following. First, we show that for classical spacing policies achieving perfect tracking requires a predictor for the states of the preceding vehicle. However, communication delays make accurate predictions of these states infeasible. To overcome this problem, we introduce delayed spacing policies that rely only on the ego-vehicle’s predicted states. First we show how actuation delays can be handled and later we show that incorporation of communication delays only requires a small extension. This novel approach simplifies the controller design process and ensures feasibility in real-world scenarios.

Secondly, we note that the use of predicted states in the spacing policy can potentially induce improper platooning behavior. To address this, we define and confine our analysis to the newly introduced class of proper spacing policies. These policies ensure proper platooning behavior and serve as a foundation for stability and performance guarantees.

Thirdly, we develop a decentralized control framework based on a predecessor–follower structure. Following the approach introduced in \cite{WijnBess20}, the follower vehicle is tasked with maintaining the desired spacing policy relative to its predecessor, independent of the predecessor’s control actions. This robust design accommodates heterogeneous platoons and ensures stability when paired with proper spacing policies. We derive necessary and sufficient conditions on spacing policies that allow for synthesizing decentralized controllers. These controllers inherently rely on predicted states, and we provide explicit designs for practical implementation.

Fourthly, we analyze three delayed spacing policies. Namely the delayed equivalents of the constant, constant headway, and extended headway spacing policies. For each policy, we identify conditions under which they qualify as proper spacing policies and ensure string stability. String stability, a critical performance criterion in vehicle platoons, guarantees that perturbations diminish as they propagate through the platoon. Our decentralized control design enables string stability through careful spacing policy design.

Finally, we validate the theoretical results through experimental tests. We present results for three controllers, each corresponding to a different amount of information available about the predecessor. These results demonstrate the practicality and robustness of our proposed approach and highlight the alignment between theoretical predictions and experimental outcomes.

The remainder of the paper is structured as follows. The motivation for predicted spacing policies and the problem formulation are provided in Section~\ref{sec_prblm}. The control design and the characterization of spacing policies for which there exists a decentralized controller is given in Section~\ref{sec_control_design}. In Section~\ref{sec_spacing} we analyze three examples of spacing policies and show that tracking of these spacing policies requires a different level of information sharing. The experimental results are provided in Section~\ref{sec_experiments}. The conclusion and recommendation for further research are given in Section~\ref{sec_conclusion}.

\section{Motivation and problem formulation}\label{sec_prblm}

Consider a platoon of $N+1$ vehicles with the dynamics
\begin{align} \label{eqn_platoondynamics}
	\begin{aligned}
		\dot{q}_i(t) &= v_i(t), \\
		\dot{v}_i(t) &= a_i(t), \\
		\tau_i\dot{a}_i(t) &= -a_i(t) + u_i(t-\phi_i),
	\end{aligned} \quad i\in\calI_N,
\end{align}
where $\calI_N:=\{0,1,..,N\}$. Here, $q_i$, $v_i$, and $a_i$ (all in $\R$) are the longitudinal position, velocity, and acceleration of vehicle $i$, respectively. The control input is given by $u_i\in\R$ and is assumed to act with an input-delay $\phi_i\in \R_+$ on the system. This input delay is required to accurately describe the response of the experimental vehicles in \cite{ploeg_2014, Lidstrom_2012, Haan_2024}. The model~(\ref{eqn_platoondynamics}) is similar to the models used in, e.g., \cite{stankovic_2000, ploeg_2014b, Haan_2024}, where the time constant $\tau_i>0$ representing the engine dynamics is not taken to be identical for each vehicle in the platoon. Thus, a heterogeneous platoon is considered. We assume that the input on the interval $[-\phi_i,0]$ is known, such that solutions for $t\geq 0$ can be computed. As a consequence, the state of the vehicle $i$ can be collected as $(x_i, u_i) \in \R^3 \times \calC^{[-\phi_i,0]}$, where $x_i = \mb q_i & v_i & a_i \mbb ^{\top}\in\R^3$. The full state of the platoon is denoted by $(x,u)$, where 
\begin{align*}
x = \mb x_0^{\top} & \cdots & x_N^{\top}\mbb^{\top}\in\R^{3(N+1)}    
\end{align*}
and 
\begin{align*}
    u=(u_0,\; \cdots,\; u_N)\in \prod_{i=0}^N \calC^{[-\phi_i,0]}.
\end{align*}

We denote the distance between the ego-vehicle $i$ and its predecessor with index $i-1$ as
\begin{align}
	\Delta_i(t) = q_{i-1}(t) - q_i(t) \label{eqn_spacing}
\end{align}
and define the desired inter-vehicle distance $\Delta_i^{\mathrm{ref}}:\R^2\rightarrow\R$ of $(v_i,a_i)$.
Such function is usually referred to as the \emph{spacing policy} and is typically formulated in terms of the velocity and acceleration of the predecessing vehicle. 
The inter-vehicle distance
$\Delta_i$ and the spacing policy $\Delta_i^\mathrm{ref}$ give rise to the spacing error $e_i$ defined~as
\begin{align}
	e_i(t) = \Delta_i(t) - \Delta_i^{\mathrm{ref}}(t).
	\label{eqn_spacingerror}
\end{align}

Although in principle a spacing policy can be defined as a general function of time, a particular subset of spacing policies is relevant in the context of platooning. This subset constitutes the set of \textit{proper} spacing policies, which is formally defined by input-to-state stability-like properties \cite{Sontag95} as follows.  
\begin{mydef}\label{def_spacing_policy}
	A function $\Delta_i^\mathrm{ref}:\R^2\rightarrow \R$ is called a \emph{proper spacing policy} for the platoon \eqref{eqn_platoondynamics}, if there exist $\beta\in \calK\calL$ and $\gamma_i, \overline \gamma_i \in \calK_\infty$ such that
   \begin{enumerate}[$i)$]
    \item $|\Delta_i(t)|\leq \beta(|\Delta_i(0)|,t)+\gamma_i(|v_{i-1}|_\infty)+\overline \gamma_i(|e_i|_\infty)$   ;
		\item$|\dot \Delta_i(t)|\leq\beta(|\dot \Delta_i(0)|,t)+ \gamma_i(|a_{i-1}|_\infty)+\overline \gamma_i(|\dot e_i|_\infty)$ ;
		\item $\,|a_{i}(t)| \,\leq \beta(\,|a_i(0)|\,,t) +\gamma_i(|a_{i-1}|_\infty)+\overline \gamma_i(|\ddot e_i|_\infty)$ .
	\end{enumerate}
\end{mydef}
The choice for proper spacing policies in the context of platooning can intuitively be motivated as follows. For a platoon driving in steady state, it is desired that each vehicle drives at the same constant velocity in the case of perfect tracking, \textit{i.e.,} $\Delta_i=\Delta_i^\mathrm{ref}$ for all $t\geq 0$. This is guaranteed in the case of proper spacing policies by part $ii)$ and $iii)$ of Definition~\ref{def_spacing_policy}. 
In fact, these conditions are slightly weaker, stating that if the spacing error remains constant, every vehicle in the platoon is driving with the same velocity in steady state.  Furthermore, for a platoon in steady-state, a bounded distance with respect to a predecessor is desired, which is ensured by part $i)$.  

Standard examples of proper spacing policies are the constant spacing policy (where $\Delta_i^{\mathrm{ref}}(t) = 0$, \textit{e.g.}, \cite{swaroop_1999}) and the constant headway policy and extended spacing policy given by
\begin{subequations}
   \begin{align}
	\Delta_i^{\mathrm{ref}}(t) &= h_v v_i(t),\label{eqn_constant_headway} \\
	\Delta_i^\mathrm{ref}(t) &= h_vv_i(t) +h_a a_i(t), \label{eqn_extended}
\end{align} 
\end{subequations}

respectively, with $h_v,h_a>0$, \textit{e.g.}, \cite{ioannou_1993,swaroop_1994,WijnBess20}.

\begin{remark}\label{rem_standstilldistance} 
    After introducing a so-called standstill distance $d_i$ in the spacing policy, the distance $q_i^\prime = q_i-d_i$ can be regarded as the deviation from the longitudinal position $q_i$, due to the linearity of the dynamics. The dynamics remain unchanged as $\dot q_i^{\prime} = v_i$  and $\Delta_i =q_{i-1}^\prime-q_i^\prime= q_{i-1}-q_i+d_i$ can be regarded as the deviation from the nominal spacing. Consequently, we can assume without loss of generality that $d_i=0$. 
\end{remark}

Due to the presence of the input delay, perfect tracking of spacing policies that depend only on current states \textit{i.e.,} $x(t)$ can be cumbersome. In the case of \textit{e.g.,} the constant headway spacing policy \eqref{eqn_constant_headway}, perfect tracking, \textit{i.e.,} $\Delta_i = \Delta_i^\mathrm{ref}$ for all $t$ implies $\ddot \Delta_i(t) = \ddot \Delta_i^\mathrm{ref}(t)$, which after a time shift $\phi_i$ is the case if and only if
\begin{align*}
	a_{i-1}(t+\phi_i)-a_i(t+\phi_i) &= h_v\dot a_i(t+\phi_i),\\
	&=\frac{h_v}{\tau_i}\bigg( u_i(t)-a_i(t+\phi_i)\bigg).
\end{align*}
Hence perfect tracking can only be achieved if future states of the preceding and following vehicle, \textit{i.e.,} the future accelerations, are known. Due to the input delay in the vehicle model~\eqref{eqn_platoondynamics} the future states of vehicle $i$ only depend on the input up to time $t$ since for $i\in\calI_N$ we have
\begin{align}\label{eqn_predictor_a}
	a_i(t+\phi_i) = e^{-\frac{\phi_i}{\tau_i}}a_i(t)+\int_{t-\phi_i}^t \frac{1}{\tau_i}e^{-\frac{t-s}{\tau_i}} u_i(s)\; ds.
\end{align}
and hence can accurately be predicted using a state-predictor. However, since the state prediction of the predecessors' acceleration or its input at time $t$ needs to be communicated in practice, an accurate estimate of $a_{i-1}(t+\phi_i)$ is in general not available at time $t$. Furthermore, if the input delay of the predecessor is smaller than the input delay of the ego vehicle, this estimate cannot be obtained at all. Hence, perfect tracking can in general not be achieved in practice.

To address this issue, we introduce linear \textit{delayed spacing policies} that incorporate actuation and communication delays. Specifically, we define spacing policies incorporating actuation delays as functional mappings $\Delta_i^\mathrm{ref}: \R^2\times\calC^{[-\phi_i,0]}\rightarrow \R$ of $(v_i,a_i)$ and $u_i$ on $[-\phi_i,0]$ of the form
\begin{align}\label{eqn_delay_spacing}
	\Delta_i^\mathrm{ref}(t) = H_ix_i(t)+ \overline H_ix_i(t+\phi_i).
\end{align}
These spacing policies depend only on current and predicted states of the follower vehicle.

In particular, we consider the \textit{delayed constant spacing policy}, the \textit{delayed constant headway policy}, and the \textit{delayed extended spacing policy}, given by \eqref{eqn_delay_spacing} where
\begin{subequations}\label{eqn_delay_spacing_pol}
	\begin{align}
        & \Delta_i^\mathrm{ref} =\int_{t}^{t+\phi_i} v_i(s)\;  ds,
        \label{eqn_dconst_sp}\\
		& \Delta_i^\mathrm{ref} = h_v v_i(t+\phi_i),\label{eqn_dconst_head_sp}\\
		& \Delta_i^\mathrm{ref} =h_v v_i(t)+h_a a_i(t+\phi_i),\label{eqn_dext_sp}
	\end{align}
\end{subequations}
respectively, for some positive constants $h_v,h_a>0$. Note $\eqref{eqn_dconst_sp}$ can also be written in the form \eqref{eqn_delay_spacing} by rewriting 
\begin{align*}
    \int_{t}^{t+\phi_i} v_i(s)\;  ds= q_{i}(t+\phi_i)-q_i(t).
\end{align*}
The spacing errors resulting from these spacing policies can be compactly written as
\begin{subequations}\label{eqn_model_spacingerror}
\begin{align}
	e_i(t) & = q_{i-1}(t)-q_i(t+\phi_i),\\
	e_i(t) &= q_{i-1}(t)-q_i(t)-h_v v_i(t+\phi_i),\\
	e_i(t) &=q_{i-1}(t)-q_{i}(t)-h_vv_i(t) -h_aa_{i}(t+\phi_i).
\end{align}
\end{subequations}
These spacing policies will shown to be proper and hence in the case of perfect tracking, \textit{i.e.}, $\Delta_i(t) =\Delta^\mathrm{ref}_i(t)$ for all $t\geq 0$, the inter-vehicle distances follow
\begin{align*}
    \Delta_i(t) = \phi_i v_i(t) \quad \text{ and } \quad \Delta_i(t) = h_v v_i(t),
\end{align*}
in steady state, \textit{i.e.,} $a_i(t)=0$, for the delayed constant headway and the delayed constant headway and delayed extended spacing policy, respectively.

Furthermore, observe that from \eqref{eqn_dext_sp}, \eqref{eqn_extended} is recovered in the case there are no delays. If in addition $h_a=0$ we obtain \eqref{eqn_constant_headway}. Due to the delayed nature of these spacing policies, they will in fact not satisfy the requirements of Definition~\ref{def_spacing_policy} for all values of $h_v$ and $h_a$. However, in the subsequent analysis characterizations in terms of $h_v$ and $h_a$ will be presented for which the spacing policies in \eqref{eqn_delay_spacing_pol} are in fact proper.

The control objectives for a controller $u_i$ of the form \eqref{eqn_controller_form} in the context of platooning can be phrased as achieving tracking and asymptotic stabilization of the spacing error while ensuring string stability. These three objectives are formally defined as follows. 
\begin{mydef}\label{def_tracking_stringstab}
	Consider the platoon (\ref{eqn_platoondynamics}) and a spacing policy $\Delta_i^{\mathrm{ref}}$. Then, a controller $u_i$ for all $i\in\calI_N$ is said to
	\begin{enumerate}[$i)$]
		\item track the spacing policy if for any $u_0(\cdot)$ and with $x(0)=0$, it holds that $e_i(t) = 0$ for all $t\geq 0$ and all $i\in\calI_N$;
		\item asymptotically stabilize the spacing policy if, for any $u_0(\cdot)$ and any $x(0)\in \R^{3(N+1)}$, it holds that
		\begin{align}
			\lim_{t\rightarrow\infty} e_i(t) = 0, \quad \forall i\in\calI_N;
		\end{align}
		\item achieve string stability if for any $u_0(\cdot)$ and with $x(0)=0$, it holds that
		\begin{align}
			\int_0^T |v_i(t)|^2 \; d t \leq \int_0^T |v_{i-1}(t)|^2 \; d t,
			\label{eqn_def_stringstab}
		\end{align}
		for all $T>0$ and all $i\in\{2,\ldots,N\}$.
	\end{enumerate}
\end{mydef}

Clearly parts $i)$ and $ii)$ of Definition~\ref{def_tracking_stringstab} correspond to the tracking performance of the platoon. An important aspect in this definition is that the properties of tracking and string stability are required to hold \emph{for any} input $u_0(\cdot)$ of the leading vehicle in the platoon. Namely, it is assumed that the lead vehicle is either manually driven or has a controller with a local control objective such as maintaining a constant speed. Thus, controllers affecting the follower vehicles (with indices $i\in\calI_N$) are responsible for satisfying the objectives in Definition~\ref{def_tracking_stringstab}. 

As mentioned in the previous section, a \emph{decentralized} approach towards controller design is pursued, in which a vehicle $i$ is responsible for achieving the desired spacing with respect to its direct predecessor $i-1$ using only local and communicated measurements. Specifically, we assume that vehicle $i$ has access to the current and past states and input of its predecessor, in addition to its own states and state predictions. Moreover, the aim is to design the control for vehicle $i$ such that it is robust with respect to the input of vehicle $i-1$. In that case, while vehicle $i$ maintains the desired inter-vehicle distance with respect to its predecessor, the latter remains free to choose its control input $u_{i-1}$ to achieve its own control objectives (e.g., tracking of a desired spacing policy with respect to vehicle $i-2$). 

Given this approach, which is similar to the control design in \cite{WijnBess20}, it is sufficient to design a controller for two consecutive vehicles. Hence, we consider a platoon consisting of only two vehicles. The state of this platoon is given by
\begin{equation}
	\begin{split}
		\xi_i &=[\begin{array}{ccc} x_{i-1}^\top &x_i^\top \end{array}]^\top,\\ 
		&= [\begin{array}{cccccc} q_{i-1} & v_{i-1} & a_{i-1} & q_i & v_i & a_i \end{array}]^\top,
	\end{split}
\end{equation}
with the corresponding dynamics
\begin{align}
	\dot{\xi}_i(t) &= \calA_i\xi_i(t) + \calB_iu_i(t-\phi_i) + \calE_iu_{i-1}(t-\phi_{i-1}),
	\label{eqn_model}
\end{align}
where
\begin{align}
	\calA_i = \left[\begin{array}{cc}  A_{i-1} & 0 \\ 0 &  A_i \end{array}\right], \;\;
	\calB_i = \left[\begin{array}{c} 0 \\  B_i \end{array}\right], \;\;
	\calE_i = \left[\begin{array}{c}  B_{i-1} \\ 0 \end{array}\right].
	\label{eqn_model_matrices}
\end{align}
Here, $ A_i$ and $ B_i$ follow directly from (\ref{eqn_platoondynamics}) and are given by
\begin{align}
	 A_i = \left[\begin{array}{ccc} 0 & 1 & 0 \\ 0 & 0 & 1 \\ 0 & 0 & -\tau_i^{-1} \end{array}\right], \qquad
	 B_i = \left[\begin{array}{c} 0 \\ 0 \\ \tau_i^{-1} \end{array}\right].
	\label{eqn_model_matricestilde}
\end{align}
Since the spacing policies and hence the spacing errors are still linear we seek a state feedback controller for vehicle $i$ of the form
\begin{align}\label{eqn_controller_form}
	u_i(t) =F_i \xi_i(t)+G_ix_i(t+\phi_i)+c_i u_{i-1}(t-\lambda) ,
\end{align}
for some $F_i\in \R^{6\times 1}$, $G_i\in \R^{3\times 1}$, $c_i,\lambda  \in \R$ and  $\lambda\geq 0$. Note that a controller of the form \eqref{eqn_controller_form} only depends on the states of vehicle $i-1$ and $i$ at time $t$, the state predictions of vehicle $i$ and a delayed input of the predecessor. This requirement on the controller implies that to compute $x(t)$ for $t\geq 0$, an initial trajectory needs to be specified, \textit{i.e.,} $x(t)$ for $t\in [-\phi_i,0]$. 

Assuming that the spacing policy takes the form \eqref{eqn_delay_spacing} the spacing error can be compactly written as 
\begin{align}\label{eqn_spacing_error_general}
	e_i(t) = \Delta_i(t) - H_i x_i(t) - \overline H_i x_i(t+\phi_i).
\end{align}
Now, a decentralized state feedback controller for vehicle $i$ given as \eqref{eqn_controller_form} leads to the closed-loop dynamics
\begin{multline}\label{eqn_closed_loop}
	\dot \xi_i(t) = \left(\calA_i+\calB_i\mb 0&G_i\mbb \right)\xi_i(t)+\calB_i F_i \xi_i(t-\phi_i)\\
	+\calB_ic_i u_{i-1}(t-\lambda-\phi_i)+\calE_iu_{i-1}(t-\phi_{i-1}),
\end{multline}
where it is noted that the control input $u_{i-1}$ is unaffected by this local controller.

 Using \eqref{eqn_closed_loop}, the objective for the design of $F_i$, $G_i$ and $c_i$ is stated in the following problem.

\begin{problem}\label{prb_statefeedback}
	Given platoon dynamics \eqref{eqn_model} and spacing error \eqref{eqn_model_spacingerror} resulting from the spacing policy $\Delta^\mathrm{ref}$, find a controller of the form \eqref{eqn_controller_form} such that the closed-loop system \eqref{eqn_closed_loop}
	satisfies the following properties for any $u_{i-1}(\cdot)$:
	\begin{enumerate}[$i)$]
		\item $\xi_i(t)=0$ on $[-\phi_i,0]$ implies $e_i(t) = 0$ for all $t\geq0$;
		\item for all initial trajectories $\xi_i\in(\calC^{[-\phi_i,0]})^6$, it holds that $\lim_{t\rightarrow\infty}e_i(t) = 0$.
	\end{enumerate}
\end{problem}
Clearly, properties $(i)$ and $(ii)$ in Problem~\ref{prb_statefeedback} correspond to the objectives of tracking and asymptotic stabilization of the spacing policy as in items $(i)$ and $(ii)$ in Definition~\ref{def_tracking_stringstab}. Although a controller that solves part $i)$ and $ii)$ in Problem~\ref{prb_statefeedback} induces zero-dynamics, \textit{i.e.,} the remaining dynamics of vehicle~$i$ when $e(t)=0$ for all $t\geq 0$, these internal dynamics are stable in the case of proper spacing policies. Indeed, after noting that the state of the platoon is given by $
    (\xi_i, u_{i-1},u_i) \in \R^6\times \calC^{[-\phi_{i-1},0]}\times \calC^{[-\phi_i,0]}
$
we can define $S: \R^6 \times \calC^{[-\phi_i,0]}\rightarrow  \R^6 \times \calC^{[-\phi_i,0]}$ by
\begin{align}\label{eqn_coordinate_transform}
	S:(\xi_i^\top ,u_{i}) \mapsto (x_{i-1}, e_i, \dot \Delta_i, a_i,u_{i-1}, u_i),
\end{align}
which constitutes a coordinate transformation, the internal dynamics of vehicle $i$ are stable if and only if the dynamics of $\dot \Delta_i$ and $a_i$ are ISS w.r.t. the behavior of vehicle $i-1$. It follows from the problem formulation and the definition of a proper spacing policy that the internal dynamics are stable. 

\begin{lemma}
	Consider the dynamics \eqref{eqn_model}. A controller that solves Problem~\ref{prb_statefeedback} for a proper spacing policy, induces input-to-state stable internal dynamics w.r.t. $a_{i-1}$. 		
\end{lemma}
\BP
Let $u_i$ be a controller that solves Problem~\ref{prb_statefeedback}. It suffices to consider the dynamics of the internal state-representation \eqref{eqn_coordinate_transform}. As a consequence of the choice of controller, $e_i(t)$ and its derivatives are stable by design and hence trivially ISS w.r.t. $a_{i-1}$. Furthermore, 
\begin{align*}
    |\dot \Delta_i(t)|&\leq \beta(|\dot \Delta_i(0)|,t)+\gamma_i(|a_{i-1}(t)|)+\bar \gamma_i(|\dot e_i|_\infty),\\
    &\leq \beta(|\Delta_i(0)|,t)+\gamma(|a_{i-1}|_\infty)+\bar \gamma_i(|\dot e_i|_\infty),
\end{align*}
proving ISS w.r.t. $a_{i-1}$ as well. The result for $a_i$ follows analogously. 
\EP 

Observe that there is no requirement on string stability in Problem~\ref{prb_statefeedback}. However, string stability is defined for a zero initial condition and hence this property depends on the zero-dynamics, \textit{i.e.,} the remaining dynamics in the case $e_i(t)=0$ for $t\geq 0$. By adopting the control objectives in Problem~\ref{prb_statefeedback}, string stability follows as a consequence of the spacing policy. As such, this property can be guaranteed by a proper design of the spacing policy.

\section{Control design and implementation}\label{sec_control_design}
In the following, we will focus on the synthesizing controllers that solve Problem~\ref{prb_statefeedback} for a particular spacing policy.

%\subsection{Control design}
Before designing controllers that solve Problem~\ref{prb_statefeedback} for the spacing policies in \eqref{eqn_delay_spacing_pol} we will present a more general result. We characterize all spacing policies for which part $i)$ of Problem~\ref{prb_statefeedback} can be solved and then show that, for these spacing policies, solvability of the second part of Problem 1 is implied by the first part. To do so, we define the relative degrees of the system \eqref{eqn_model} with respect to an output of the form \eqref{eqn_spacing_error_general}.

\begin{mydef}
	The smallest integers $\rho_i$ and $\overline \rho_i$ such that
	\begin{align*}
		H_i A_i^{\rho_i-1}B_i \neq 0, \qquad \overline H_i A_i^{\overline \rho_i -1}B_i\neq 0.
	\end{align*}
	are called \emph{the relative degrees} of the system \eqref{eqn_platoondynamics} with respect to the output \eqref{eqn_spacingerror}.
\end{mydef}
\begin{theorem}\label{thrm_existence_sp}
	Consider the dynamics \eqref{eqn_model} with a spacing policy of the form \eqref{eqn_delay_spacing}. Then part $i)$ of Problem~\ref{prb_statefeedback} has a solution if and only if the relative degrees $\rho_i$ and $\overline \rho_i$ satisfy $\overline \rho_i< \rho_i$, or $\overline \rho_i =3$ implies $H_ix_i = -q_i$.
\end{theorem}
\BP
$(\Rightarrow)$ Let $u_i$ be the controller that solves Problem~\ref{prb_statefeedback}. Then if $e_i(t)=0$ for all $t\geq 0$ also, $\dot e_i(t) =\ddot e_i(t)=\dddot{e_i}(t)=0$ for all $t\geq 0$. Then $\dot e_i(t) =0 $ implies
\begin{multline}\label{eqn_derror_perfect}
	v_{i-1}(t)-v_i(t) = H_i A_ix_i(t) +H_iB_i u_i(t-\phi_i)\\+\overline H_i A_i  x_i(t+\phi_i)+\overline H_i B_i u_i(t).
\end{multline}
In the case $\rho_i=1$ and $\overline \rho_i >1$, it follows directly from \eqref{eqn_derror_perfect} that
\begin{multline*}
    u_i(t-\phi_i) = (H_iB_i)^{-1}\Big(v_{i-1}(t)-v_i(t)-H_iA_i x_i(t)\\-\overline H_iA_ix_i(t+\phi_i) \Big).
\end{multline*}
Hence, after shifting the controller with $\phi_i$ time, the controller depends on the predicted velocity of vehicle $i-1$, which contradicts the assumption that $u_i$ takes the form \eqref{eqn_controller_form}. Hence we conclude $\overline \rho_i\leq \rho_i$. 

In the case $\overline \rho_i =1$, substituting the controller form \eqref{eqn_controller_form} for $u_i(t-\phi_i)$ yields
\begin{multline*}
     u_i(t) = (\overline H_iB_i)^{-1} \Big(v_{i-1}(t)-v_i(t)
     \\-H_i(A_i-B_iG_i)x_i(t)   -\overline H_i A_i x_i(t+\phi_i) \\- H_i B_i (F_i \xi_i(t-\phi_i) -c_iu_{i-1}(t-\lambda-\phi_i)\big).
\end{multline*}
Due to the controller form \eqref{eqn_controller_form}, the controller is assumed to be independent of the states $\xi_i(t-\phi_i)$ and thus it follows that $F_i=0$. Then substituting the controller form for $u_i(t)$ with $F_i=0$ in \eqref{eqn_derror_perfect} yields 
\begin{multline*}
    u_{i}(t-\phi_i) =(H_iB_i)^{-1}\Big(v_{i-1}(t)-v_i(t)-H_iA_i x_i(t)\\\qquad \quad \;\qquad \quad \;\qquad \quad\;\quad-(\overline H_iA_i+G_i)x_i(t+\phi_i)\\ -c_iu_{i-1}(t-\lambda)\Big),
\end{multline*}
which again shows that after a shift $\phi_i$, the input depends on the predicted velocity on the velocity $v_{i-1}(t+\phi_i)$ of the predecessor, which again violats the assumption that $u_i$ is of the form \eqref{eqn_controller_form}. Consequently, $\rho_i >1$. 

In the case $\overline \rho_i=1$ the statement is proven and hence we continue with the assumption that $\overline \rho_i >1$. Observe that $\ddot e_i(t)=0$ for all $t\geq 0$ implies
 \begin{multline}\label{eqn_dderror_perfect}
 	a_{i-1}(t)-a_i(t) = H_i A^2_ix_i(t) +H_iA_iB_i u_i(t-\phi_i)\\+\overline H_i A^2_i  x_i(t+\phi_i)+\overline H_iA_i B_i u_i(t)
 \end{multline}
Using similar arguments as in the previous case, it follows that $\rho_i>2$. Hence it remains to show the case where $\overline \rho_i=3$. In that case $\dddot e_i(t)=0$ for all $t\geq 0$ it follows that
 \begin{align}\label{eqn_ddderror_perfect}
 \begin{aligned}
	\dot a_{i-1}(t)-\dot a_i(t)& =\frac{1}{\tau_{i-1}}\left (u_{i-1}(t-\phi_{i-1})-a_{i-1}(t)\right)\\
	&\quad \qquad \quad \quad-\frac{1}{\tau_i}(u_{i}(t-\phi_i)-a_i(t)),\\
	&= H_i A^3_ix_i(t) +H_iA_i^2B_i u_i(t-\phi_i)\\
	&\quad +\overline H_i A^3_i  x_i(t+\phi_i)+\overline H_iA_i^2 B_i u_i(t).
    \end{aligned}
\end{align}
Suppose $\rho_i =3$. Then the controller is independent of past states of the predecessor if and only if $H_iA_i^2B_i =-\frac{1}{\tau_i}$. However, this implies $H_i = \mb -1 &0&0 \mbb$, which is equivalent to $H_ix_i = -q_i$. 

$(\Leftarrow)$ An input that ensures $e(t)=0$ for all $t\geq 0$ under the given assumptions follows trivially from \eqref{eqn_derror_perfect} in the case of $1= \overline \rho_i <\rho_i$ and \eqref{eqn_dderror_perfect} when $2= \overline \rho_i <\rho_i$ and from \eqref{eqn_ddderror_perfect} if $\overline \rho_i =3$ and $H_i x_i = -q_i$. 
\EP

The conditions in Theorem~\ref{thrm_existence_sp} allow for insightful interpretation. First of all, the result implies that $\overline H_i\ne 0$ is a necessary condition for a spacing policy to admit a solution to Problem~\ref{prb_statefeedback}. Stated differently, predicted states must be used in the spacing policy to enable perfect tracking. This observation partially explains the problems experienced in mitigating the effect of the delays when $\overline H_i=0$. For spacing errors independent of these predicted states may still be stabilizable, but the input $u_{i-1}$ can not be decoupled from the spacing error by a feedback of the form \eqref{eqn_controller_form}. Secondly, the relative degrees must satisfy $\rho_i>\overline \rho_i$. This means that the input $u_i(t)$ must appear sooner, \textit{i.e.,} in a lower order derivative of the spacing error, than the delayed input $u_i(t-\phi_i)$.

The following result states that the conditions of Theorem~\ref{thrm_existence_sp}
are sufficient to also solve part $ii)$ of Problem~\ref{prb_statefeedback}.
\begin{theorem}\label{thrm_existence_stab}
	Consider the dynamics \eqref{eqn_model} with a spacing policy of the form \eqref{eqn_delay_spacing}. If part $i)$ of Problem~\ref{prb_statefeedback} can be solved, then part $ii)$ can be solved simultaneously and the controllers are given by
    \begin{align}\label{eqn_controller_rho1}
		&	\begin{aligned}
			u_i(t) =& -\left (\overline H_i B_i \right)^{-1}\Big(v_{i-1}(t)-v_{i}(t)-H_iA_i x_i(t)\\
			&  -\overline H_i A_i x_i(t+\phi_i)+k_p e_i(t)\Big)
		\end{aligned}
        \end{align}
and
\begin{align}\label{eqn_controller_rho2}
		&\begin{aligned}
			u_i(t) =& -\left (\overline H_iA_i B_i \right)^{-1}\Big(a_{i-1}(t)-a_{i}(t)-H_iA_i^2 x_i(t)\\
			&  -\overline H_i A_i^2 x_i(t+\phi_i)+k_p e_i(t)+k_d\dot e_i(t)\Big),
		\end{aligned}
	\end{align}
	for $\overline \rho_i =1$ and $\overline \rho_i=2$, respectively, and 
	\begin{align}\label{eqn_controller_rho3}
    \begin{aligned}
		u_i(t) =& -\left (\overline H_iA_i B_i \right)^{-1}\bigg(-\overline H_i A_i^2 x_i(t+\phi_i)\\
        &\quad \qquad  \frac{1}{\tau_{i-1}}\Big(u_{i-1}(t-\phi_{i-1})-a_{i-1}(t)\Big)\\
		& \quad \qquad \qquad  +k_p e_i(t)+k_d\dot e_i(t)+k_{dd} \ddot e_i(t)\bigg),
        \end{aligned}
	\end{align}
	for $\overline \rho_i =3$. Furthermore, the gains satisfy $k_p>0$, $k_d>0$, $k_{dd}>0$ and $0<k_{p}k_d-k_{dd}$.
\end{theorem}
\BP
The controllers that solve Problem~\ref{prb_statefeedback} ensure $e_i(t)=0$ for all $t\geq T$ if $e_i(T)=0$ and hence the derivatives of $e_i(t)$ are zero for $t\geq T$. Consequently, these controllers can be regarded as input-output linearizing controllers. By definition, the spacing error is of the form \eqref{eqn_spacing_error_general} and hence it is, as well as its derivatives, a linear combination of $\xi_i(t)$ and $x_i(t+\phi_i)$. A such, the controllers \eqref{eqn_controller_rho1}, \eqref{eqn_controller_rho2} and \eqref{eqn_controller_rho3} are of the form \eqref{eqn_controller_form}. Substituting the controllers yields that the error dynamics are given by 
\begin{align*}
	\dot e_i(t) &=-k_p e_i(t),\\
	\ddot e_i(t)&=-k_d\dot e_i(t)-k_p e_i(t),\\
	\dddot e_i(t)&= -k_{dd} \ddot e_i(t)-k_de_i(t)-k_pe_i(t).
\end{align*}
Consequently, these controllers ensure $e_i(T)=0$ for all $T\geq t$ when $e_i(t)=0$. The conditions on the gains follow straightforwardly from the Routh-Hurtwitz criterion for stability of these error dynamics.
\EP 

\begin{remark}
 Although the results of Theorem~\ref{thrm_existence_sp} and \ref{thrm_existence_stab} are formulated within a linear framework, these results can be straightforwardly generalized to nonlinear spacing policies and nonlinear error dynamics. To do so, nonlinear techniques from \textit{e.g.} \cite{book_nijmeijer_1990,khalil2015nonlinear} are to be exploited. However, as the proofs for nonlinear results are more complex and would detract from the main focus of this paper, they are considered beyond its scope. 
\end{remark}

\begin{remark}
In the case of both the delayed constant and
delayed constant headway spacing policies, the controller depends on information of the predecessor, which must be communicated. However, due to communication delays, the controller for the following vehicle can only access delayed states of the predecessor at the current time $t$. To address this issue, an alternative definition of the inter-vehicle distance such as:
\begin{align}\label{eqn_alt_spacing}
   \overline  \Delta_i(t) =q_{i-1}(t-\varphi_i)-q_i(t),
\end{align}
could be considered, where $\varphi_i>0$ denotes the communication delay between vehicles $i-1$ and $i$. A spacing error involving this alternative inter-vehicle distance \eqref{eqn_alt_spacing} can be induced by adopting the spacing $\overline \Delta_i(t)$ and a spacing policy of the form 
\begin{align}\label{eqn_spacing_pol_alt}
    \overline \Delta_i^\mathrm{ref}(t)=\int_{t-\varphi_i}^t v_{i-1}(s) \; ds + \Delta_i^\mathrm{ref}(t),
\end{align} 
with $\Delta_i^\mathrm{ref}$ of the form \eqref{eqn_delay_spacing}. This adapted spacing policy is proper if $\Delta_i^\mathrm{ref}$ is proper and a controller that solves Problem~\ref{prb_statefeedback} can be designed along similar lines as in Theorem~\ref{thrm_existence_stab} and only requires delayed states of vehicle $i-1$. However, including this communication delay in the spacing policy leads to a larger inter-vehicle distance in steady state, since
\begin{align*}
    \overline \Delta_i (t)&= \Delta_i(t) -\varphi_i v_{i-1}(t)= \Delta_i(t) -\varphi_i v_{i}(t).
\end{align*}

In practice only the acceleration $a_{i-1}(t)$ is required to be communicated, while other states such as \textit{e.g.,} $q_{i-1}(t)-q_i(t)$ and $v_{i-1}(t)-v_i(t)$ are measured on board of vehicle $i$ and do not suffer from communication delays. As such, dealing with communication delays requires a more detailed analysis and is therefore considered outside of the scope. Furthermore, in our experimental setup, the communication delay is an order of magnitude smaller than the actuator delay and the spacing $\Delta_i$ and velocity difference $\dot \Delta_i$ are measured using radar. Since radar measurements are assumed not to be subject to communication delays, we chose not to adopt this delayed spacing policy in our analysis.
\end{remark}

\section{Analysis of spacing policies}\label{sec_spacing}
In the following, we will apply the results obtained in the previous section to the spacing policies in \eqref{eqn_delay_spacing_pol}. In particular, we will show that for these spacing policies the conditions of Theorem~\ref{thrm_existence_sp} are satisfied. The controllers that solve Problem~\ref{prb_statefeedback} follow form Theorem~\ref{thrm_existence_stab} and hence it remains to show that the spacing policies are in fact proper spacing policies.

\subsection{Delayed constant spacing policy}
We will start with the analysis of the delayed constant spacing policy \eqref{eqn_dconst_sp}. A controller that achieves tracking and asymptotic stabilization follows straightforwardly from Theorem~\ref{thrm_existence_stab} as follows. 
\begin{corollary}\label{cor_constant}
	Consider the platoon \eqref{eqn_platoondynamics} and the delayed constant spacing policy \eqref{eqn_dconst_sp}. A controller of the form \eqref{eqn_controller_form} solves part $i)$ of Problem~\ref{prb_statefeedback} if and only if
\begin{multline}\label{eqn_controller_constant}
		u_i(t)=\frac{\tau_i}{\tau_{i-1}}\bigg(-a_{i-1}(t)+u_{i-1}(t-\phi_{i-1})\bigg)+a_i(t+\phi_i)\\+\tau_i\bigg(k_pe_i(t)+k_d \dot e_i(t) +k_{dd} \ddot e_i(t)\bigg).
	\end{multline}
	for some constants $k_p>0$, $k_d>0$, $k_{dd}>0$ and $0<k_pk_d-k_{dd}$.
\end{corollary}

The result of Corollary~\ref{cor_constant} shows that the delayed constant spacing policy can only be tracked and stabilized if the delay $\phi_{i-1}$ and the engine constant $\tau_{i-1}$ are known to the predecessor unless a heterogeneous platoon is considered. 

It remains to show that the delayed constant spacing policy is a proper spacing policy. However, for this spacing policy, this is trivially the case. 
\begin{theorem}
    Consider the delayed constant spacing policy \eqref{eqn_dconst_sp}. The spacing policy is a proper spacing policy.
\end{theorem}
\BP 
After writing 
\begin{align*}
    \Delta_i(t) &= \int_{t-\phi_i}^t v_{i-1}(s)\; ds+e_i(t-\phi_i),\\
    \dot \Delta_i(t) &= \int_{t-\phi_i}^t a_{i-1}(s)\; ds+\dot e_i(t-\phi_i),\\
    a_i(t) &= a_{i-1} (t-\phi_i)+\ddot e_i(t-\phi_i).
\end{align*}
satisfaction of the requirements in Definition~\ref{def_spacing_policy} follow. Hence, the spacing policy is proper. 
\EP 

\begin{theorem}
    Consider the platoon \eqref{eqn_platoondynamics} with the delayed constant headway spacing policy. Let the inputs $u_i$ be given by \eqref{eqn_input_ch} for every $i\in \calI_N$. Then the platoon is string stable.
\end{theorem}
\BP
As a consequence of the choice of controller, for a zero initial trajectory, we have $v_{i}(t+\phi_i)=v_{i-1}(t)$ for all $t$. As such, in the Laplace domain, we obtain
\begin{align*}
    \hat v_i (s) = e^{-\lambda s} \hat v_{i-1}(s)=T(s)\hat v_{i-1}(s).
\end{align*}
Consequently, since $|e^{-\lambda i\omega}|=1$, we can conclude, using the equivalence of the $H_\infty$ norm \cite{ploeg_2014,book_dullerud_2000}, that 
\begin{align*}
    |T_i(i\omega)|_\infty:=\sup_{\omega\in \R^+} |T_i(i\omega)| =\sup_{v_{i-1}\neq 0} \frac{|v_i|_{\calL_2}}{|v_{i-1}|_{\calL_2}}\leq 1.
\end{align*}
This proves the desired result. 
\EP 

\subsection{Delayed constant headway}
Next we analyze the delayed constant headway spacing policy \eqref{eqn_dconst_head_sp}. Similarly to the delayed constant spacing policy, the controller design follows straightforwardly from Theorem~\ref{thrm_existence_stab}, leading to the following result.

\begin{corollary}\label{cor_constant_h}
	Consider the platoon \eqref{eqn_platoondynamics} and the delayed constant headway spacing policy \eqref{eqn_dconst_head_sp}. A controller of the form \eqref{eqn_controller_form} solves Problem~\ref{prb_statefeedback} if and only if
	\begin{multline}\label{eqn_input_ch}
		u_i(t) =  a_i(t+\phi_i)\\
		+\frac{\tau_i}{h_v}\bigg(a_{i-1}(t)-a_i(t)+k_p e_i(t)+k_d\dot e_i(t)\bigg).
	\end{multline}
	for some $k_p>0$ and $k_d>0$. 
\end{corollary}

In order to synthesize the controller given in Theorem~\ref{cor_constant_h}, the acceleration of the predecessor is required to be available to the follower. This state is typically measured on board and communicated V2V. As such, this controller can be regarded as a CACC controller. 

It remains to characterize when the delayed constant headway spacing policy is a proper spacing policy. Whereas the delayed constant spacing policy was trivially proper, this is not the case for the constant headway-spacing policy. 

\begin{theorem}\label{thrm_dch_proper}
	Consider the delayed constant headway spacing policy \eqref{eqn_dconst_head_sp}. The spacing policy is a proper spacing policy if and only if $2\phi_i<h_v\pi$.
\end{theorem}
\BP 
The proof is provided in the Appendix. 
\EP

Next we show when string stability is achieved. In this paper string stability is defined for a zero initial condition and the controller is required to achieve tracking of the spacing policy. As such, due to the control design, the achievement of string stability does not depend on the controller gains but on the choice of spacing policy. This is reflected in the following characterization.

\begin{theorem}\label{thrm_dch_stringstab}
	Consider the platoon \eqref{eqn_platoondynamics} with the delayed constant headway spacing policy. Let the inputs $u_i$ be given by \eqref{eqn_input_ch} for every $i\in \calI_N$. Then the platoon is string stable if and only if $h_v\geq 2\phi_i$ 
\end{theorem}
\BP 
The proof is provided in the Appendix. 
\EP

\subsection{Delayed extended headway spacing policy}

Finally, consider the extended headway spacing policy. A controller that solves Problem~\ref{prb_statefeedback} is characterized in the next theorem. 

\begin{theorem}
	Consider the platoon \eqref{eqn_platoondynamics} and the delayed extended spacing policy \eqref{eqn_dconst_head_sp}. A controller of the form \eqref{eqn_controller_form} solves Problem~\ref{prb_statefeedback} if and only if
	\begin{multline}\label{eqn_controller_extended}
		u_i(t)=\frac{\tau_i}{h_a} \bigg( v_{i-1}(t)-v_i(t)-h_va_i(t)+k_pe_i(t)\bigg)\\+   a_i(t+\phi_i),
	\end{multline}
	for some $k_p>0$. 
\end{theorem}

Note that a controller that tracks the delayed extended spacing policy can be interpreted as an ACC controller, as it only requires knowledge of the states of the follower vehicle and relative states. As such, all information required can be measured using onboard measurement equipment. This is in contrast to the delayed constant and delayed constant headway spacing policy, which requires wireless communication for the synthesis of this controller.  

\begin{remark}
    By choosing the gain $k_p=\tau_i^{-1}$ the controller \eqref{eqn_controller_extended} reduces to
    \begin{multline*}
        u_i(t) = \frac{\tau_i}{h_a} \bigg(  v_{i-1}(t)-v_i(t)-h_va_i(t) \\
        +q_{i-1}(t)-q_i(t) -h_v v_i(t)\bigg).
    \end{multline*}
    Consequently, for this choice of $k_p$, the controller does not depend on predicted states. This means that in principle, no knowledge of the actuation delay is required and no state predictor. 
\end{remark}

It remains to characterize when the delayed extended headway spacing policy is a proper spacing policy. Similar to the delayed constant headway spacing policy, only certain constants $h_v$ and $h_a$ will induce a proper spacing policy. 

\begin{theorem}\label{thrm_dex_proper}
	Consider the extended delayed spacing policy \eqref{eqn_dext_sp}. The spacing policy is a proper spacing policy if and only if for some $\omega \in \left(0, \frac{\pi}{2}\right)$ the headways $h_v$ and $h_a$ satisfy
	\begin{align}\label{eqn_stability_region}
		\frac{\phi_i h_v}{h_a} &< \omega \sin(\omega),\\
        \frac{\phi_i^2}{h_a} &< \omega^2 \cos(\omega).
	\end{align}
\end{theorem}
\BP 
The proof is provided in the Appendix. 
\EP

\begin{figure}[ht]
	\centering
	\includegraphics[width=\linewidth]{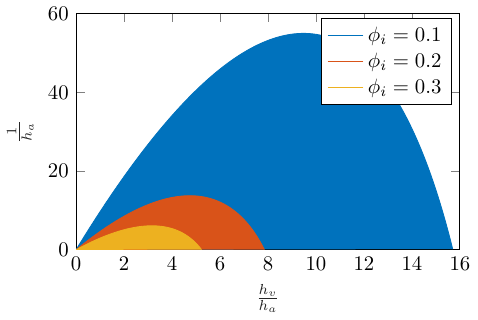}
	\caption{The region where the values of $\frac{h_v}{h_a}$ and $\frac{1}{h_a}$ induce properness of the spacing policy \eqref{eqn_dext_sp} for different delays $\phi_i$. }
	\label{fig_stab_region}
\end{figure}

The region for which the delayed extended headway spacing policy is a proper spacing policy, shrinks as the the delay increases. To visualize this, Figure~\ref{fig_stab_region} shows the stability region for various choices of $\phi_i$. 
\begin{remark}
    The result of Theorem~\ref{thrm_dex_proper} implies that
\begin{align*}
    2\phi_ih_v < h_a \pi,
\end{align*}
which is not surprising as in a steady state, \textit{i.e.,} $a_i(t)=0$ the spacing policy coincides with the constant headway spacing policy. Furthermore, note that string stability implies
\begin{align*}
    \frac{\phi_i^2}{h_a} < \frac{\pi^2}{4},
\end{align*}
indicating that there exists an infimum for $h_a$.  
\end{remark}

The following result provides sufficient conditions for the extended delayed spacing policy to induce string stability given controllers solving Problem~\ref{prb_statefeedback}.
\begin{theorem}\label{thrm_dex_stringstab}
    Consider the platoon \eqref{eqn_platoondynamics} with the delayed extended headway spacing policy. Let the inputs $u_i$ be given by \eqref{eqn_input_ch} for every $i\in \calI_N$. Then the platoon is string stable if
    \begin{align*}
	h_a \geq   2h_v\phi_i, \quad \text{ and } \quad h_v^2\geq 2h_a.
\end{align*}
\end{theorem}
\BP 
The proof is provided in the Appendix. 
\EP
\begin{remark}
Observe that setting $k_p=\tau_i^{-1}$ in controller \eqref{eqn_controller_extended} effectively reduces the controller to a PD-controller for the constant headway spacing policy without considering actuation delays. Remarkably, Theorem~\ref{thrm_dex_stringstab} demonstrates that such a PD-controller, despite not accounting for actuation delays explicitly, can still achieve string stability. 
\end{remark}
\begin{remark}
    The result of Theorem~\ref{thrm_dex_stringstab} shows that the condition $h_v^2\geq 2h_a$ is necessary for string stability, but not sufficient. However, it is also a sufficient condition if $\phi_i=0$, \textit{c.f.} Theorem~17 in \cite{WijnBess20}.
\end{remark}

\section{Experimental validation}\label{sec_experiments}
\begin{figure}
\center
\includegraphics[width = 0.75\linewidth]{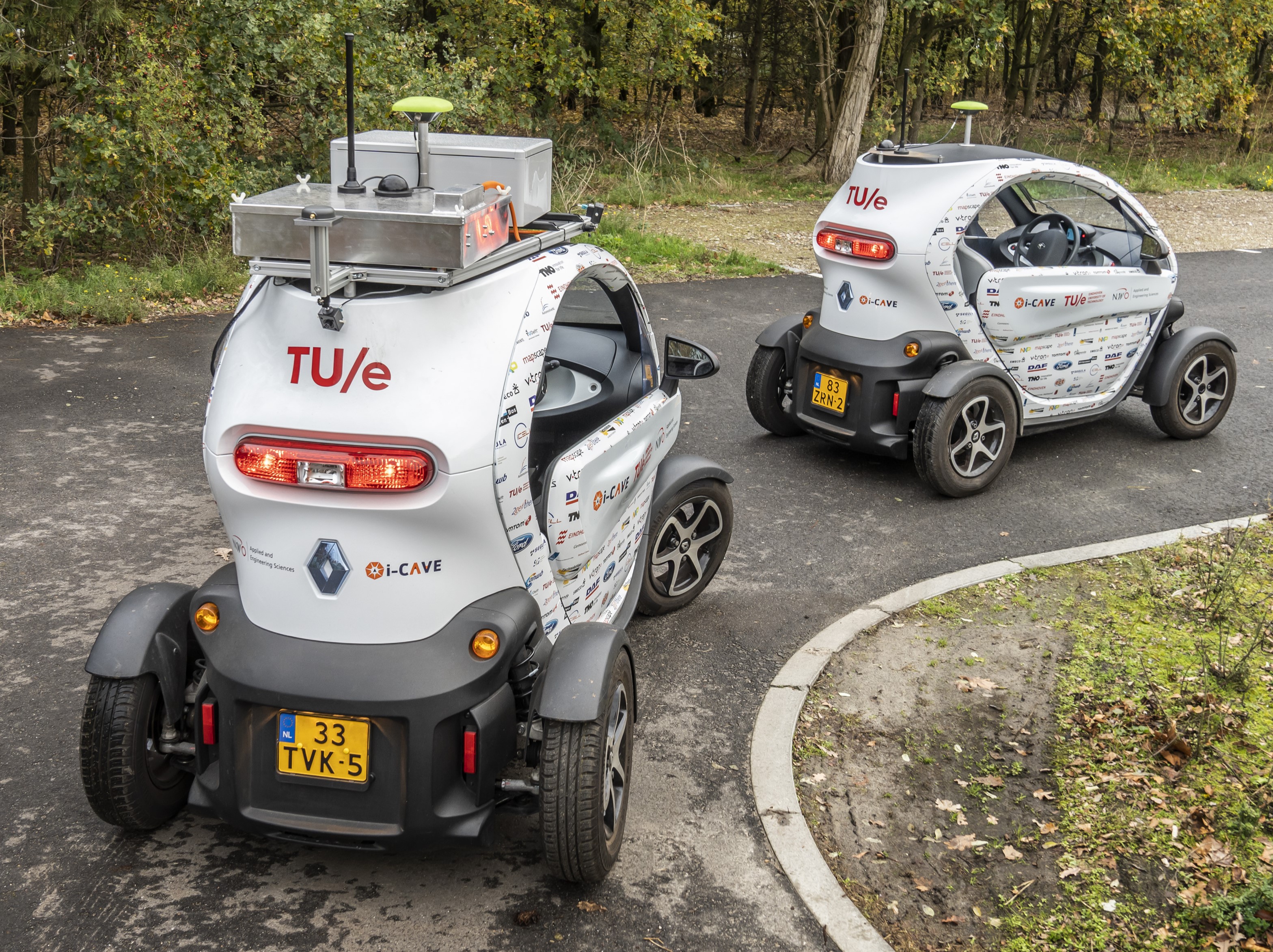}
\caption{A platoon consisting of two test vehicles.}
\label{fig:experimental-vehicles}
\end{figure}

To evaluate the proposed framework and spacing policies introduced in Section~\ref{sec_spacing}, experiments are performed with a platoon consisting of two full-scale vehicles as shown in \figref{fig:experimental-vehicles}. These experimental vehicles are two-seated electric vehicles that have been modified so that the resulting platform exhibits drive-by-wire technology. The throttle and braking commands, determined by the controller running on a real-time platform, are executed by the vehicle. In this article, only the general architecture and characteristics related to the longitudinal control of the vehicles is discussed. The lateral control during the experiments is performed by a human safety driver. More details regarding the automation of the vehicles can be found in \cite{Hoogeboom_2020}.

\begin{figure}
\center
\includegraphics[width = 0.75\linewidth]{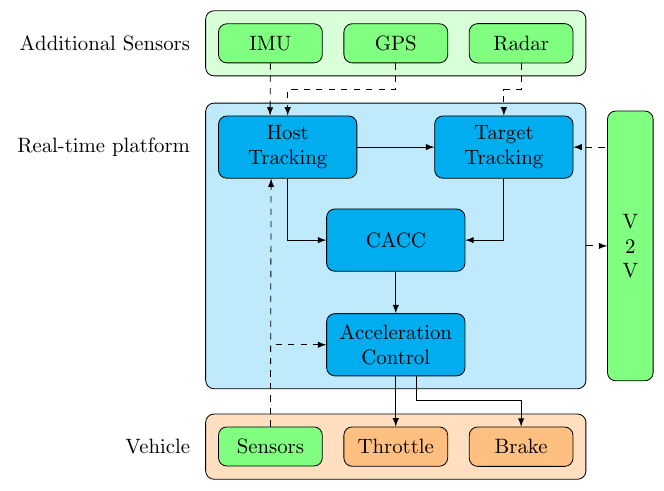}
\caption{Schematic overview of the experimental vehicle's control system and sensors.}
\label{fig:experimental-schematic}
\end{figure}

\subsection{Experimental setup}
The schematic control structure of the experimental vehicles is depicted in Figure~\ref{fig:experimental-schematic}. The control system that operates at a base frequency of 100 Hz has access to the (original) CAN network of the vehicle, which provides (in the scope of longitudinal control) measurements of the rear-axle (rotational) velocity. In addition, the vehicle is equipped with an IMU for acceleration measurements. A forward facing radar measures the relative distance and relative velocity of objects, denoted by range and range rate respectively, at a rate of 16.7 Hz. By utilizing a V2V router, the vehicles are able to broadcast information to surrounding vehicles using the ITS G5 standard~\cite{Strom_2011}. However, this wireless communication link introduces communication latency. In the experiments, the V2V communication was configured to transmit at a rate of 25 Hz, resulting in an average communication latency of 20 ms. Additionally, the vehicles are equipped with a GPS that provides position and groundspeed measurements, as well as a common clock among vehicles that is used to synchronize the local logs. The sensors that are used to obtain the required measurements for the controllers are listed in \tabref{tab:sensors-information}. 

\renewcommand{\arraystretch}{1.2}
\begin{table}[htb]
\centering
\caption{Measurements performed by associated sensors on the ego vehicle.}
\label{tab:sensors-information}
\begin{tabular}{| l | l l |}\hline 
Sensor & \multicolumn{2}{l| }{Measurement} \\ \hline \hline
\multirow{2}{*}{Radar} & Range & $\Delta_i $ \\
 & Range rate & $\dot{\Delta}_i$  \\ \hline
Speedometer & Ego velocity & $v_i$ \\ \hline
IMU & Ego acceleration & $a_i$ \\ \hline
\multirow{3}{*}{V2V} & Leader input & $u_{i-1}$ \\
& Leader velocity & $v_{i-1}$ \\
& Leader acceleration & $a_{i-1}$\\ \hline
\end{tabular}
%\vspace{-1em}
\end{table}
 
\begin{figure}
\resizebox{.95\linewidth}{!}{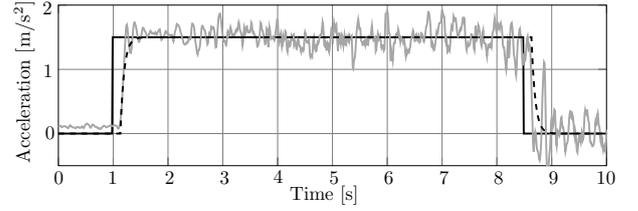}
\caption{Measured acceleration response (gray) of experimental vehicle subject to a step input $u_i(t)$ (solid black) compared to the model \eqref{eqn_platoondynamics} response with $\tau_i = 0.067~\text{s}, \phi_i = 0.15~\text{s}$ (dashed black), obtained from \cite{Hoogeboom_2020}.}
\label{fig:measured-acceleration-response}
\end{figure}

\paragraph{Longitudinal control}
The longitudinal model \eqref{eqn_platoondynamics} that is adopted in the controller design, is the result of a feedback linearization, taking into account the aerodynamic drag, viscous friction and rolling resistance. Using this feedback linearizing controller, given an acceleration setpoint $u_i(t)$, the driving torque is obtained that is required to achieve this acceleration. An inverse motor mapping and (configurable) brake distribution are used to transform this driving torque into a setpoint for the electric motor of the vehicle. We chose to perform all braking during the experiments on the electric motor, to have an identical positive and negative acceleration response. The resulting vehicle response as depicted in \figref{fig:measured-acceleration-response} can be modeled by \eqref{eqn_platoondynamics} where the vehicle's parameters are $\tau_i = 0.067~\text{s}$ and $\phi_i = 0.15~\text{s}$. Note that this model only holds in a limited acceleration range. However, limit situations (which require excessive braking) are considered outside the scope of our controllers and should be handled by systems designed specifically for that particular task.

\paragraph{Predictor implementation}
Each of the controllers requires a prediction of the ego vehicle's states. Since the real-time control system operates in discrete time, the control input $u_i(t)$ is applied at discrete times to the vehicle. Assuming a zero-order hold on the input
\begin{equation} \label{eq:zoh}
u_i(t) = u_i(kT_s) \text{ for } kT_s \leq t < kT_s + T_s, \ k\in \mathbb{N},    
\end{equation}
and assuming the delay $\phi_i = d T_s$ for some integer\footnote{For readability an integer multiple of the sampling time is assumed for the delay. For non-integer multiples a similar approach can be followed, by splitting the integral associated with the input.} $d$, we can obtain an exact expression of a dynamic system 
\begin{equation*}
    \dot{x}(t) = Ax(t) + Bu(t-\phi),
\end{equation*}
at the sampling instances by 
\begin{equation*}
x(kT_s + T_s) = \Phi x(kT_s) + \Gamma u(kT_s - dT_s),
\end{equation*}
where $\Phi = e^{AT_s}$, $\Gamma = \int_{0}^{T_s} e^{A\sigma} \mathrm{d}\sigma B$. We use these discrete time dynamics to obtain a prediction of the future states as
\begin{equation}\label{eq:predictor-discrete}
    \hat{x}(kT_s + dTs) = \Phi ^d x(kT_s) + \sum_{j=1}^d \Phi^{j-1} \Gamma u(kT_s - jT_s).
\end{equation}

For the ego vehicle, this results in an exact prediction, as we know the sampling time and input delay to the system. 

\begin{remark}
The implementation of the controller \eqref{eqn_controller_constant} requires absolute position measurements and the associated prediction $q_i(t + \phi_i)$ to compute the error. Since such measurements are generally not available in practice, we chose to define the error in a frame relative to vehicle $i$ in the real-time implementation of the controller.
\end{remark}

\subsection{Experiment design}
To perform the experiments, we initialize the platoon in an equilibrium position at standstill (i.e., zero velocity) on a straight road. The leader vehicle is controlled by a cruise controller which aims to control the vehicle to a reference velocity $v_\text{ref}$. Once the reference velocity is reached in a steady state, the cruise control is disabled and positive and negative acceleration pulses are prescribed (open-loop) to the leader vehicle. The follower vehicle is controlled using the controllers discussed in Section~\ref{sec_spacing}. Each experiment is performed with the controller tuning and headways that theoretically result in (string) stability as listed in \tabref{tab:controller-tuning}.

\renewcommand{\arraystretch}{1.2}
\begin{table}[htb]
\centering
\caption{Controller tuning and headways used in the experiments.}
\label{tab:controller-tuning}
\begin{tabular}{| l | l |}\hline 
Controller & Parameter choice \\ \hline \hline
\multirow{3}{*}{Constant spacing \eqref{eqn_controller_constant} }& $k_p = \tfrac{1}{\tau_i}$ \\
& $k_d = \tfrac{3}{\tau_i} $ \\
& $k_{dd} = \tfrac{3}{\tau_i} $ \\ \hline
\multirow{3}{*}{Constant headway \eqref{eqn_input_ch} } & $h_v = 0.4 $\\
& $k_p = 0.2$ \\ 
& $k_d = 0.7 - \tau_i k_p$ \\ \hline
\multirow{3}{*}{Extended headway \eqref{eqn_controller_extended} } & $h_v = 1.2 $ \\
& $h_a = 0.25 $ \\
& $k_p = 0.2$ \\ \hline
\end{tabular}
\end{table}

\subsection{Experimental results}
\begin{figure}
\begin{center}
\begin{subfigure}{0.99\linewidth}
\includegraphics[width=\linewidth]{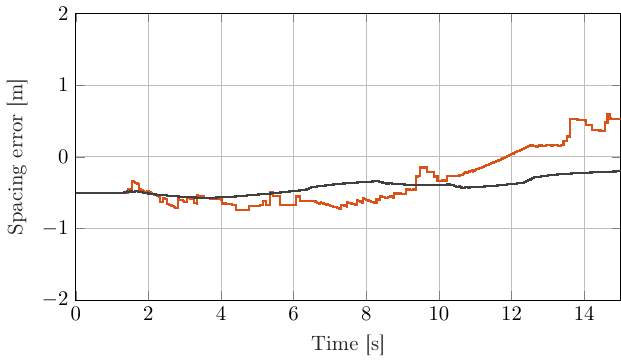}
\caption{Error $e_i(t)$.}
\end{subfigure}
\begin{subfigure}{0.99\linewidth}
\includegraphics[width=\linewidth]{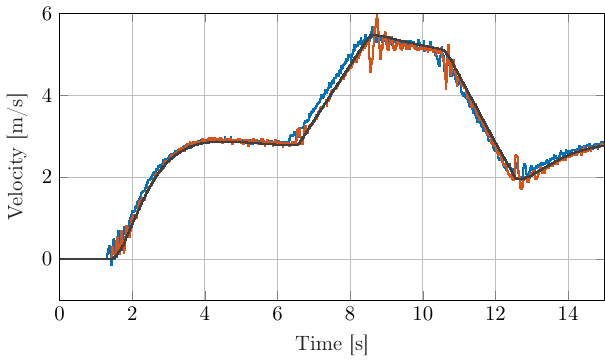}
\caption{Velocity response.}
\end{subfigure}
\begin{subfigure}{0.99\linewidth}
\includegraphics[width=\linewidth]{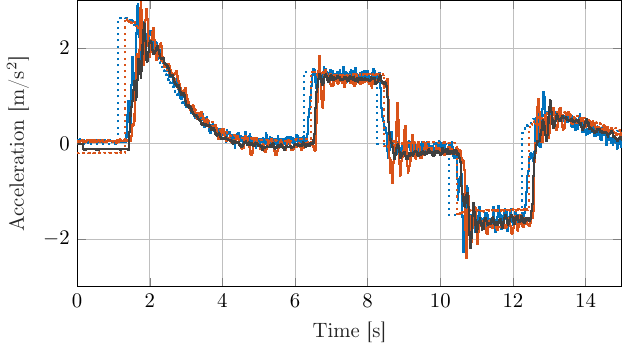}
\caption{Acceleration response (solid) and setpoint (dashed).}
\end{subfigure} 
\caption{Measured experimental response (\protect\legendFollower), and simulated response (\protect\legendFollowerSimulated) of vehicle deploying controller \eqref{eqn_controller_constant} for tracking of the delayed constant spacing policy, with tuning from \tabref{tab:controller-tuning}, to follow leader vehicle (\protect\legendLeader).} 
\label{fig:experimental-response-constant}
\end{center}
\end{figure}

\begin{figure}
\begin{center}
\begin{subfigure}{0.99\linewidth}
\includegraphics[width=\linewidth]{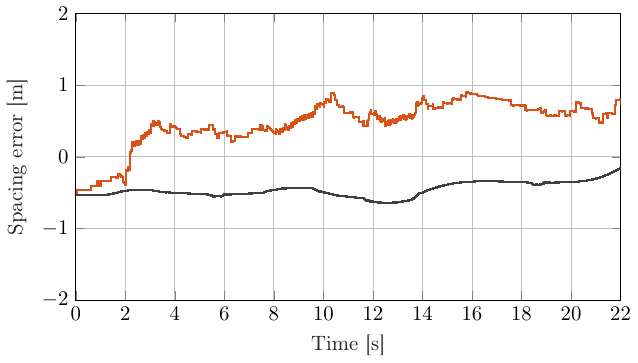}
\caption{Error $e_i(t)$.}
\end{subfigure}
\begin{subfigure}{0.99\linewidth}
\includegraphics[width=\linewidth]{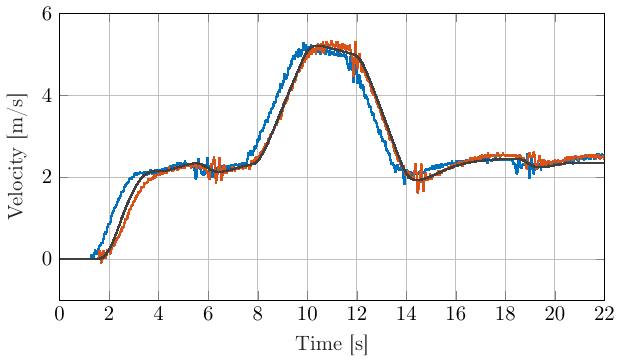}
\caption{Velocity response.}
\end{subfigure}
\begin{subfigure}{0.99\linewidth}
\includegraphics[width=\linewidth]{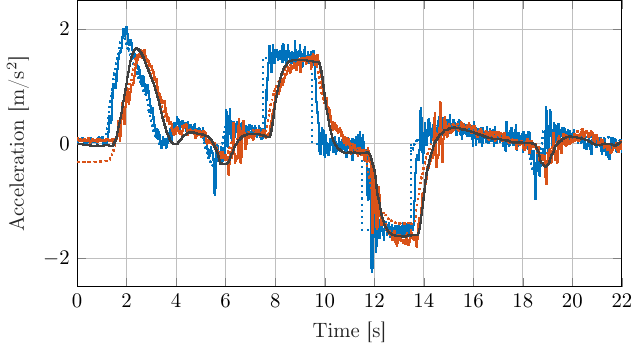}
\caption{Acceleration response (solid) and setpoint (dashed).}
\end{subfigure} 
\caption{Measured experimental response (\protect\legendFollower), and simulated response (\protect\legendFollowerSimulated) of vehicle deploying controller \eqref{eqn_input_ch} for tracking of the delayed constant headway spacing policy, with tuning from \tabref{tab:controller-tuning}, to follow leader vehicle (\protect\legendLeader).} 
\label{fig:experimental-response-headway}
\end{center}
\end{figure}

\begin{figure}
\begin{center}
\begin{subfigure}{0.99\linewidth}
\includegraphics[width=\linewidth]{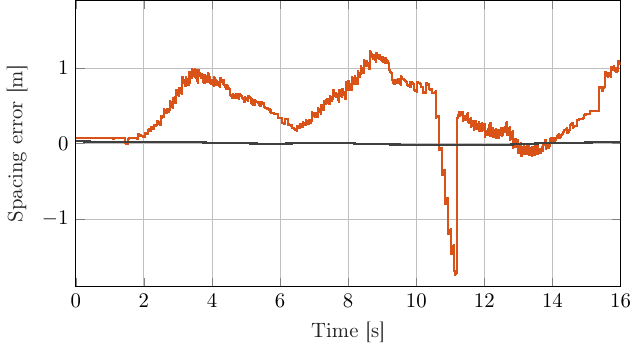}
\caption{Error $e_i(t)$.}
\end{subfigure}
\begin{subfigure}{0.99\linewidth}
\includegraphics[width=\linewidth]{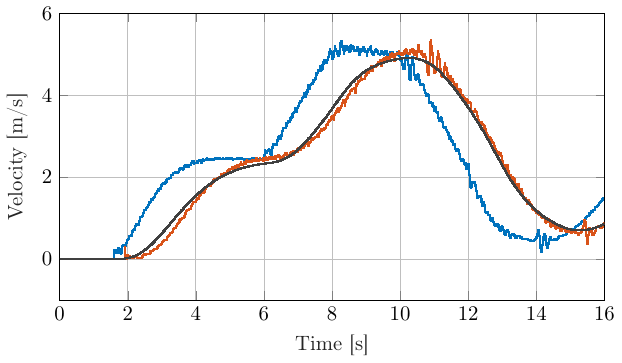}
\caption{Velocity response.}
\end{subfigure}
\begin{subfigure}{0.99\linewidth}
\includegraphics[width=\linewidth]{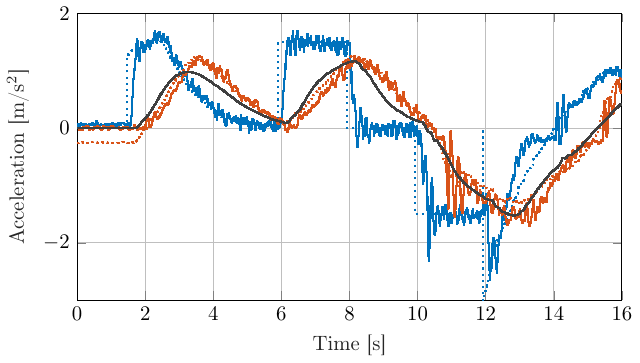}
\caption{Acceleration response (solid) and setpoint (dashed).}
\end{subfigure} 
\caption{Measured experimental response (\protect\legendFollower), and simulated response (\protect\legendFollowerSimulated) of vehicle deploying controller \eqref{eqn_controller_extended} for tracking of the delayed extended spacing policy, with tuning from \tabref{tab:controller-tuning}, to follow leader vehicle (\protect\legendLeader).} 
\label{fig:experimental-response-extended}
\end{center}
\end{figure}

The measured experimental results for the delayed constant spacing controller \eqref{eqn_controller_constant}, the delayed constant headway \eqref{eqn_input_ch} and the delayed extended spacing policy \eqref{eqn_controller_extended} are shown in \figref{fig:experimental-response-constant},~\figref{fig:experimental-response-headway} and~\figref{fig:experimental-response-extended} respectively. 

In the acceleration response of all three experiments, a negative acceleration setpoint can be observed at the start of the experiment. This negative setpoint is due to the initial distance at which the vehicles were initialized, which was slightly smaller than the desired standstill distance, resulting in a small error. However, the experimental vehicles were physically not able to drive backward to correct for this error. Instead, the negative acceleration setpoint results in a braking torque which causes the vehicle to remain stationary until the leader vehicle starts to drive. 

Due to the relatively small gains that we chose on the errors (and their derivatives), the feedback controller is not able to achieve zero errors, but they remain small over the course of the experiment. These small feedback gains highlight the effectiveness of our proposed framework, as, despite the minimal effort of the feedback controller, the errors remain small. Additionally, we clearly see disturbance decoupling of the error with respect to the input $u_{i-1}$ of the leader vehicle for both the delayed constant spacing \eqref{eqn_input_ch} and delayed constant headway \eqref{eqn_controller_constant} controllers. For the extended spacing controller \eqref{eqn_controller_extended} however, this seems not to be the case. We suspect this is due to the $v_{i-1} - v_i$ term in the input, which is obtained through the radar measurement. The radar uses an object tracking algorithm, which introduces additional dynamics that we did not account for. However, despite this lesser performance with respect to the error, the response of the controller does exhibit string stable behavior. 
Specifically, both from the response from standstill, and the pulses while driving, all controllers exhibit string stable behavior in the measured response.

Additionally, we compare the experimental results to a simulation of a vehicle modeled according to \eqref{eqn_platoondynamics} employing the respective controllers, which act on the measured experimental states of the leader vehicle. The results show how well the measured experimental results match the modeled response. The largest mismatch between measured response and simulated response can be observed in the implementation of the extended delayed spacing controller \eqref{eqn_controller_extended} as depicted in \figref{fig:experimental-response-extended}. As this controller uses only relative measurements that are performed with the radar on board of vehicle $i$, e.g., it measures the velocity difference $v_{i-1} - v_i$. In the simulation, we did not account for any dynamics that such automotive radar potentially induces. Consequently, this could explain the difference between simulation and measurements. Nevertheless, the errors remain small and the response is string stable.

Finally, the fact that all three spacing policies and the associated controllers transition to experimental responses, confirming the theoretical results, demonstrates the power of the delayed spacing policy approach.

\section{Conclusion and recommendations}\label{sec_conclusion}
This paper investigated delay-based spacing policies for vehicle platoons, which account for actuation delays. The concept of proper spacing policies was introduced, defining spacing policies that are relevant in the context of platooning. A decentralized control approach was presented, requiring each vehicle’s control objective to be robust to the control input of its predecessor. Using input-output linearization techniques, necessary and sufficient conditions were derived for the existence of controllers that achieve both tracking and asymptotic stabilization of the desired spacing policy. Three delayed equivalents of commonly used spacing policies were analyzed and shown to be proper. Furthermore, conditions were formulated under which these spacing policies ensure string stability for tracking controllers. The proposed methods were validated experimentally using full-size vehicles, demonstrating that both tracking of the spacing policy and string stability can be achieved in practice.

The presented framework could also accommodate communication delays with some generalization of the results. While these extensions are relatively straightforward, they were not explored in this paper to maintain clarity and due to space limitations. Future work will focus on integrating delayed measurements into the framework and analyzing their impact.

Furthermore, in this paper, the controllers have been tested only at relatively low velocities. Their performance at speeds typical of real traffic conditions remains to be evaluated. Another important research direction is to investigate whether incorporating additional states of the preceding vehicle into the spacing policy—such as in \eqref{eqn_delay_spacing}, can enhance string stability in practice. Additionally, engine constants, as well as actuation delays, can vary over time and may depend on vehicle velocity. Since the current framework does not account for these uncertainties, future work should explore methods to incorporate them.

Furthermore, how to integrate performance metrics beyond tracking and string stability remains an open question. One of the overarching goals of platooning is to reduce fuel consumption, which is currently only implicitly addressed by assuming that smaller inter-vehicle distances reduce drag. However, tracking the spacing policy may not be optimal for fuel efficiency, due to the transient behavior. Therefore, an important future direction is to develop tracking controllers that ensure string stability while minimizing fuel consumption and other relevant metrics. 

\section*{Appendix}

In this appendix, we provide the proofs of the theorems in Section~\ref{sec_spacing}. The proofs rely mainly on theory for delay differential equations. To briefly introduce the results required to prove our statements, consider the delay differential equation
\begin{align} \label{eqn_delay_sys}
    \dot x_i(t) =\overline  A_0 x(t)+\overline A_1 x(t-\phi)+Bu_i.
\end{align}
Globally asymptotic stability of \eqref{eqn_delay_sys} is characterized as follows \cite{michiels_2011}. 
\begin{proposition}\label{prop_delay_stab}
    Consider the linear time-delay system \eqref{eqn_delay_sys} and assume $u_i\equiv 0$. The solution $x(t)\equiv 0$ is asymptotically stable if and only if the characteristic equation $p_\chi(\lambda)$ given by
    \begin{align*}
       p_\chi(\lambda)= \det\left (\lambda I - \overline A_0-\overline A_1e^{-\lambda \phi}\right),
    \end{align*}
    has all its roots in the open left half plane. 
\end{proposition}
Input to state stability of such a delayed differential equation is characterized as follows \cite{pepe2006lyapunov,chaillet2023iss}.
\begin{theorem}\label{thrm_delay_iss}
The linear time-delay system \eqref{eqn_delay_sys} is ISS with respect to $u_i$ if and only if the trivial solution $x(t)\equiv 0$ of the unforced system (\textit{i.e.,} $u_i\equiv 0$) is
asymptotically stable.
\end{theorem}

Given these results, we now return to proving the results in Section~\ref{sec_spacing}.

\noindent \textbf{Proof of Theorem~\ref{thrm_dch_proper}.} Observe that given the spacing policy \eqref{eqn_dconst_head_sp} we obtain after rearranging terms in the spacing error and its derivative
\begin{align*}
	\Delta_i(t) &= e_i(t)+h_v v_i(t+\phi_i),\\
	&= e_i(t)-h_v \dot \Delta_i(t+\phi_i)+h_v v_{i-1}(t+\phi_i),\\
	\dot \Delta_i(t) &= \dot e_i(t) +h_v a_i(t+\phi_i).
\end{align*}
Since $e_i(t)$ and $\dot e_i(t)$ are globally asymptotically stable by design, they are trivially ISS w.r.t. $v_{i-1}$. Hence we can conclude that $\Delta_i$ is ISS with w.r.t. $\dot \Delta_i$ and $v_{i-1}$ and that $\dot \Delta_i(t)$ in turn is ISS with respect to $a_{i}$. Hence it suffices to show that $a_i$ is ISS with respect to $a_{i-1}$. To do so, consider the closed loop dynamics resulting from the input \eqref{eqn_input_ch}:
\begin{align*}
	\ddot e_i(t) &= -k_p e_i(t)-k_d\dot e_i(t)\\
	\dot a_i(t) &=-\frac{1}{h_v}\bigg (-a_i(t-\phi_i)+a_{i-1}(t-\phi_i) \\
	&\qquad +k_pe_i(t-\phi_i)+k_d \dot e_i(t-\phi_i)\bigg).
\end{align*}
After defining $y = \mb e_i & \dot e_i & a_i\mbb^\top$, we can rewrite these dynamics as
\begin{align*}
	\dot y(t) =\overline A_0 y(t)+\overline A_1y(t-\phi_i)+\overline B a_{i-1}(t-\phi_{i})
\end{align*}
where 
\begin{align*}
	\overline A_0 = \mb 0&1&0\\-k_p&-k_d&0\\0&0&0\mbb , \qquad \overline A_1 = \mb 0&0&0\\0&0&0\\ \frac{k_p}{h_v}&\frac{k_d}{h_v}&-\frac{1}{h_v} \mbb.
\end{align*}
These dynamics are ISS with respect to $a_{i-1}$ if and only if the trajectory $x(t)\equiv0$ is globally asymptotically stable. This is the case if and only if the characteristic equation
\begin{multline*}
	p_\chi(\lambda)= \det(\lambda I-A_0-A_1 e^{-\phi_i\lambda}),\\
	= ( \lambda +h_v^{-1}e^{-\phi_i \lambda})(\lambda^2+k_d\lambda +k_p),
\end{multline*}
has all its zeros in the open left half plane, \textit{i.e.,} $\C^-$.

Suppose $\lambda = \rho+\omega  i $ is a root of $ e^{-\phi_i \lambda}+h_v\lambda$. Then
\begin{align*}
	e^{-\rho\phi_i} (\cos(\omega\phi_i)-i\sin(\omega\phi_i)= - h_v (\rho+\omega i).
\end{align*}
As $p_\chi(\lambda)$ is a complex analytic function, $p_\chi(\lambda)=0$ implies $p_\chi(\bar \lambda)=0$ an hence we can assume without loss of generality that $\omega>0$. Suppose $\rho \geq 0$. Then $\cos(\omega\phi_i)\leq 0<\sin(\omega\phi_i)$. Hence we conclude that
\begin{align*}
	\omega \phi_i\in S :\bigcup_{k=0}^\infty \left \{\left [\frac{k\pi}{2\phi_i},\frac{k\pi}{\phi_i}\right )\right \}.
\end{align*}
Since the real part and the complex part sum up to zero, it follows that
\begin{align*}
	f(\omega):=\frac{\sin(\omega\phi_i)}{\omega}=h_v e^{\rho\phi_i}.
\end{align*}
Observe that $\frac{d}{d\omega} f(\omega)<0$ and $f\left(\frac{\pi}{2\phi_i}\right )=\frac{2\phi_i}{\pi }$ and thus it follows that
\begin{align*}
	h_v \leq h_v e^{\rho\phi} =\omega^{-1} \sin(\omega\phi_i) \leq \frac{2\phi}{\pi}.
\end{align*}
Hence it follows that $\frac{h_v}{\phi_i}\leq \frac{2}{\pi}$ if $\rho\geq 0$. Therefore, we can conclude that $\rho <0$ if and only if
\begin{align*}
	\frac{\phi_i}{h_v}< \frac{\pi}{2},
\end{align*}
which proves the desired result. 
\EP 

\noindent \textbf{Proof of Theorem~\ref{thrm_dch_stringstab}.}

Due to the choice of controller, $e_i(0)=0$ implies $e_i(t)=0$ for all $t\geq 0$. Consequently, for an initial trajectory, the spacing error is initialized and the dynamics are given by
\begin{align*}
	h_v a_i(t+\phi_i) = v_{i-1}(t)-v_i(t),
\end{align*}
Taking the Laplace transform results in the transfer function given by
\begin{align}
	\frac{\hat v_i(s)}{\hat v_{i-1}(s)}= T(s) = \frac{1}{e^{s\phi_i}h_v s+1}.
\end{align}
Consequently, $|T(s)|_\infty\leq 1$ if and only if
\begin{align*}
	\sup_{\omega\in \R^+} |i\omega h_v + e^{-i\phi_i\omega}| \geq 1.
\end{align*}
A direct computation shows
\begin{align*}
	| i \omega h_v +e^{-i\omega\phi}|^2 & = (\omega h_v -\sin(\omega\phi_i))^2+\cos(\omega\phi_i)^2,\\
	&= \omega^2 h_v^2-2\omega h_v \sin(\omega\phi_i) +1.
\end{align*}
Hence $|T(i\omega)|_\infty\leq 1$ if and only if for all $\omega>0$
\begin{align*}
	\omega^2 h_v^2 \geq 2 \omega h_v\sin(\omega \phi_i),
\end{align*}
or alternatively written, if and only if 
\begin{align}\label{eqn_condition_stringstab}
	\frac{h_v}{2\phi_i} \geq \frac{\sin(\omega\phi_i)}{\omega \phi_i}.
\end{align}
Clearly, $h_v\geq 2\phi_i$ implies satisfaction of the inequality \eqref{eqn_condition_stringstab} for all $\omega >0$, therefore showing sufficiency of the condition. However, since the 
\begin{align*}
	\sup_{\omega\in \R^+}\frac{\sin(\omega \phi_i)}{\omega\phi_i}=1
\end{align*}
and the inequality holds for all $\omega>0$ in the case of a string stable platoon, $h_v\geq 2\phi_i$ is also necessary. 
\EP

\noindent \textbf{Proof of Theorem~\ref{thrm_dex_proper}.} 
Observe that given the spacing policy \eqref{eqn_dext_sp} we obtain after rearranging terms in the spacing error and its derivative
\begin{align*}
	\Delta_i(t) &= e_i(t)+h_v v_i(t)+h_a a_i(t+\phi_i),\\
	&= e_i(t)-h_v \dot \Delta_i(t+\phi_i)+h_v v_{i-1}(t)+h_aa_i(t+\phi_i),\\
	\dot \Delta_i(t) &= \dot e_i(t) +h_v a_i(t)+h_a\dot a_i(t+\phi_i).
\end{align*}
Consequently, we can conclude ISS of $\Delta_i$ w.r.t. $e_i$, $\dot \Delta_i$ and $v_{i-1}$ and $a_i$. Similarly, ISS of $\dot \Delta_i$ w.r.t. $a_i$ and $\dot a_i$ can be concluded. Hence properness of the spacing policy is proven if ISS of $a_i$ and $\dot a_i$ w.r.t. $a_{i-1}$ can be concluded. 

To do so, observe that it follows from the definition of the spacing error that
\begin{align}\label{eqn_intern_diff}
	h_a \ddot a_i(t+\phi_i) = a_{i-1}(t)-a_i(t)+h_v\dot a_i(t)-\ddot e_i(t).
\end{align}
Since $\dot e_i =-k_p e_i$ we can conclude that $\ddot e_i$ is globally asymptotically stable and hence inherently ISS with respect to $a_{i-1}$.  Due to the assumption that $u_i(t)$ solves Problem~\ref{prb_statefeedback} it follows from Theorem~\ref{thrm_delay_iss} that $a_i$ and $\dot a_i$ are input-to-state stable with respect to $a_{i-1}$ if and only if 
\begin{align*}
	h_a \ddot a_i(t+\phi_i)=-a_i(t)+h_v\dot a_i(t)
\end{align*}
yields a globally asymptotically stable delayed differential equation. This is the case if and only if its characteristic equation defined by
\begin{align*}
	p_\chi(\lambda)= h_a e^{\phi_i\lambda}\lambda^2+h_v \lambda +1=0,
\end{align*}
has all its roots $\lambda\in \C^-$. 
Let $\lambda = \rho+\omega i$. Then using Eulers formula $e^{-\lambda \phi_i} = e^{-\rho\phi_i}(\cos(\omega\phi_i)-i\sin(\omega \phi_i))$, a direct computation shows that $\lambda $ is a root of $ p_\chi(\lambda)$ if and only if $h_v$ and $h_a$ satisfy
\begin{align}
	\begin{aligned}\label{eqn_boundary}
		\frac{h_v}{h_a}&=   -\frac{1}{\omega}e^{\rho \phi_i}\bigg((\rho^2-\omega^2)\sin(\omega \phi_i)+2\rho \omega\cos(\omega \phi_i)\bigg),\\
		\frac{1}{h_a}& =-e^{\rho \phi_i}\bigg((\rho^2-\omega^2)\cos(\omega \phi_i)-2\rho \omega\sin(\omega \phi_i)\bigg)\\
		&\qquad-\frac{h_v}{h_a} \rho.
	\end{aligned}
\end{align}
Consequently, \eqref{eqn_intern_diff} is a globally asymptotically stable delay differential equation and $a_{i}$ is input-to-state stable with respect to $a_{i-1}$ if and only if $\frac{h_v}{h_a}$ and $\frac{1}{h_a}$ are confined to the stability region where the boundary is given by \eqref{eqn_boundary} with $\rho=0$, \textit{i.e,} the curve in the $\left (\frac{h_v}{h_a}, \frac{1}{h_a} \right)$ plane defined by
\begin{align*}
	\left(\frac{h_v}{h_a}, \frac{1}{h_a}\right)=  \left\{\big(\omega\sin(\omega \phi_i), \omega^2 \cos(\omega \phi_i)\big) \mid \omega \in \bigg(0,\frac{\pi}{2\phi_i}\bigg) \right\}.
\end{align*}
This is the case if and only if for some $ \omega \in \left(0,\frac{\pi}{2}\right)$, 
	\begin{align*}
		\frac{\phi_i h_v}{h_a} <  \omega \sin( \omega), \quad \text{and} \quad 
        \frac{\phi_i^2}{h_a} <  \omega^2 \cos( \omega).
	\end{align*}
This proves the desired result. 
\EP

\noindent \textbf{Proof of Theorem~\ref{thrm_dex_stringstab}.} 
Due to the choice of controller, $e_i(0)=0$ implies $e_i(t)=0$ for all $t\geq 0$. Consequently, for an initial trajectory, the spacing error is initialized and the dynamics are given by
\begin{align*}
    h_a a_i(t+\phi_i) = -h_v v_i(t)-q_i(t)+q_{i-1}(t).
\end{align*}
The transfer function can be computed by taking the Laplace transform and is given by
\begin{align*}
	T(s)&= \frac{1}{h_ae^{\phi_i \lambda}\lambda^2+h_v\lambda +1}.
\end{align*}
A direct computation shows that $\sup_{\omega\in \R^+} |T(i\omega)|\geq 1$ if and only if
\begin{align*}
	\sup_{\omega\in \R^+}|-h_a\omega^2 e^{i\omega \phi_i} +h_v\omega i +1| \geq 1.
\end{align*}
After substituting Eulers formula for complex numbers, we obtain
\begin{multline*}
	-h_a\omega^2 e^{i\omega \phi_i}+h_v\omega i +1 \\=(1 -h_a\omega^2 \cos(\omega \phi_i))+i(h_v\omega-h_a\omega^2 \sin(\omega\phi_i))
\end{multline*}
Consequently $\sup_{\omega\in \R^+}|T(i\omega)|\leq 1$, if
\begin{align*}
	h_a^2 \omega^4+h_v^2 \omega^2 \geq 2h_a \omega^2 \cos(\omega\phi_i)+2h_ah_v \omega^3 \sin(\omega\phi_i),
\end{align*}
for all $\omega>0$. This is the case if and only if
\begin{align*}
	h_a^2 \omega^2 +h_v^2\geq 2h_a \cos(\omega\phi_i) +2h_ah_v\omega \sin(\omega\phi_i).
\end{align*}
Using the equivalence of the $H_\infty$ norm yields the desired sufficient condition.

Over approximating terms by $\cos(\omega\phi_i)\leq 1$ and similarly $\sin(\omega\phi_i)\leq \omega \phi_i$ yields
\begin{align*}
	h_a^2 \omega^2 +h_v^2\geq 2h_a+2h_ah_v \phi_i \omega^2\geq 0,
\end{align*}
or equivalently
\begin{align*}
	(h_a^2-2h_vh_a\phi_i )\omega^2 +h_v^2\geq 2h_a.
\end{align*}
Consequently, the inequality holds for all $\omega$ if
\begin{align*}
	h_a \geq   2h_v\phi_i, \quad \text{ and } \quad h_v^2 \geq 2h_a.
\end{align*}
This proves the result. 
\EP

\bibliographystyle{acm}
\bibliography{sys_con,bibliography}

\end{document}